\documentclass[12pt]{llncs}

\usepackage{amssymb}
\usepackage{amsmath}
\usepackage{amsfonts}

\newcommand{\Z}{\mathbb{Z}}
\newcommand{\R}{\mathbb{R}}

\newcommand{\Q}{\mathbb{Q}}
\newcommand{\N}{\mathbb{N}}
\newcommand{\F}{\mathbb{F}}

\usepackage{mathtools}
\DeclarePairedDelimiter{\floor}{\lfloor}{\rfloor}

\DeclarePairedDelimiter\abs{\lvert}{\rvert}

\usepackage{tikz}
\usepackage{tikz-3dplot}
\usetikzlibrary{patterns}
\usepackage{tikz-cd}
\usepackage{adjustbox}

\usepackage{tabulary}
\usepackage{algorithm}
\usepackage{algpseudocode}

\algnewcommand\And{\textbf{and}}
\algnewcommand\Or{\textbf{or}}

\usepackage{vmargin}
\setpapersize{A4}
\setmarginsrb{2.5cm}{2.0cm}{2.5cm}{2.0cm}{0pt}{1.0cm}{0pt}{1.0cm}

\title{A New Angle on Lattice Sieving for the Number Field Sieve}

\author{Gary McGuire, Ois\'in Robinson}

\institute{UCD School of Mathematics and Statistics\\ University College Dublin\\ Ireland}

\tikzcdset{every label/.append style = {font = \large}}

\usepackage{etoolbox}
\patchcmd{\abstract}{\small}{}{}{}

\begin{document}
\maketitle

\begin{abstract}
Lattice sieving in two or more dimensions has proven to be an indispensable practical aid
in integer factorization and discrete log computations involving the number field sieve.
The main contribution of this article is to show that a different method of lattice enumeration 
in three dimensions will provide a significant speedup.
We use the  successive
minima and shortest vectors  of the lattice instead of transition vectors to iterate through lattice points.
We showcase the new method by a record computation in a 133-bit subgroup of $\F_{p^6}$,
with $p^6$ having 423 bits. Our overall timing nearly $3$ times faster than the previous record of 
a 132-bit subgroup in a 422-bit field.
The approach generalizes to dimensions 4 or more, overcoming a key obstruction
to the implementation of the tower number field sieve.  
\end{abstract}

\section{Introduction}
The most widely adopted public-key cryptography algorithms in current use are critically
dependent on the (assumed) intractability of either the integer factorization problem 
(IFP), the finite field discrete logarithm problem (DLP) or the elliptic curve discrete
logarithm problem (ECDLP).  The most effective known attacks against IFP and DLP use
the same basic algorithm, namely the Number Field Sieve (NFS).  This algorithm has
subexponential complexity in the input size.  On the other hand, all known methods to
attack the ECDLP in the general case have exponential complexity.  However there are
special instances of the ECDLP which can be attacked by effectively transferring the
problem to a finite field, allowing the NFS to be used.  For example such instances 
arise in the context of pairing-based cryptography, where certain elliptic curves can
be used to realize `Identity-Based Encryption' (IBE).  There is a trade-off between the
reduced security due to the size of the finite field on which the security is dependent,
and increased efficiency of the pairing arithmetic.  The optimal parameters have been the
subject of intense scrutiny over the last few years, which have seen a succession of
improvements in the NFS for the DLP in the medium characteristic case.  This is directly
relevant in the case of pairings, where the finite field on which the security of the
protocol depends is typically a small degree extension of a prime field.

A key part of the NFS is lattice sieving.
The main contribution of this article is to demonstrate that different methods of lattice enumeration
can make a significant difference to the speed of lattice sieving.

This paper is organized as follows.
In section 2, we give a very brief overview of the Number 
Field Sieve algorithm in the medium-characteristic
case.  
A  more detailed explanation can be found in \cite{MR3775580}.
One of the main bottlenecks of this algorithm is lattice sieving, which involves enumerating points
in a (low-dimensional) lattice.  
We  propose in Section 3 a straightforward idea to significantly increase enumeration speed in dimensions 3 and above.
The idea is to change the angle of planes that are sieved through in order to reduce the 
number of planes. This idea has been used before for lattice enumeration in a sphere, however 
it has not been applied successfully to lattice sieving for the NFS.  We show that the idea 
can work well by using integer linear programming to find an initial point for iteration in 
a plane within the sieve cuboid. 
In Section 4 we propose
a novel method to amortize memory communication overhead which applies regardless of dimension.
In section 5 we give details of a new record discrete log computation in $\F_{p^6}$.
The previous record due to Gr\'{e}my et al  \cite{MR3775580} had $p^6$ with 422 bits, and this paper
has $p^6$ with 423 bits. We deliberately chose a field
size just one bit larger because this allows a direct comparison of methods and timings.
In Section 6 we present a record pairing break with the same prime $p$.
Finally we conclude in section 7 and mention some possible future research ideas.

\section{Number Field Sieve}
\noindent We start by describing NFS in the most naive form suitable for computing discrete logs in
$\F_{p^n}$.  Consider the following commutative diagram:\\
\begin{figure}[H]
\adjustbox{scale=0.8,center}{%
	\begin{tikzcd}[row sep=small, column sep=tiny, style={font=\large}]
	 & a-bx \in \Q[x] \arrow{dl}[swap]{x \mapsto \alpha} \arrow{dr}{x \mapsto \beta} & \\
	[+45pt] \Q[x]/\langle{f_0(x)\rangle} \cong \Q(\alpha) \arrow{dr}[swap]{\alpha \mapsto m} & & \Q[x]/\langle{f_1(x)\rangle} \cong \Q(\beta) \arrow{dl} \\
	[+45pt] & (\Z/p\Z)[x]/\langle{\psi(x)\rangle} \cong \F_{p^n}
	\end{tikzcd}
}
\caption{Commutative diagram of NFS for discrete log in $\F_{p^n}$.}
\label{fig:nfscommute}
\end{figure}
\noindent The polynomials  $f_0(x)$ and $f_1(x)$
are irreducible in $\Z [x]$ of degree $n$, and they  define the number fields $\Q(\alpha)$ and $\Q(\beta)$ respectively.
We require that  $f_0(x)$ and $f_1(x)$, when reduced modulo $p$,
share a factor $\psi(x)$ of degree $n$ which is irreducible over $\F_p$.  This defines
the finite field $\F_{p^n}$ as $ (\Z/p\Z)[x]/\langle{\psi(x)\rangle}$.
Usually $\psi(x)$ is simply the reduction of $f_0(x)$ modulo $p$.

For a bound $E$, we inspect many pairs of integers $(a,b)$ with $0 < a \leq E$ 
and $-E \leq b \leq E$ in the hope of finding
many pairs such that
\[
\text{Res}\left(f_0, a - bx\right) \hspace{5mm} \text{and} \hspace{5mm} \text{Res}\left(f_1, a - bx\right)
\]
are both divisible only by primes up to a bound $B$.  
In 3-dimensional sieving, the pairs $(a,b)$ corresponding to $a-bx$ become
triples $(a,b,c)$ corresponding to $a+bx+cx^2$.

The NFS has four main stages - polynomial selection, sieving, linear algebra, descent.  For further details
see \cite{MR3344923},  \cite{MR3775580}, \cite{MR2422170}, \cite{GremyPhD}.  
Recently, new variations of NFS have been described where the
norms (i.e. resultants) are even smaller in certain fields, see \cite{MR3565295}, \cite{MR3598073}.

\subsection{Lattice Sieving}
\noindent 

The `special-$\mathfrak{q}$' lattice sieve, originally due to J.M. Pollard \cite{MR1321220} is outlined first.
Let $q$ be a rational prime, let $r$ be an integer with $f_0(r) \equiv 0 \mod{q}$,
and let $\mathfrak{q} = \langle q,\theta-r \rangle$ be
an ideal of $K = \Q(\theta) \cong \Q[x]/\langle f_0 \rangle$ lying over $q$.
We look for (integral) ideals of $K$ that are divisible by $\mathfrak{q}$, and we do this by
looking for ideals whose  norm  is divisible by $q$. 
We also would like the norm to be divisible by many other small primes $p$.
We fix $q$ and iterate over all $p$ in the factor base using a sieve.

We do this in three dimensions as follows.
We  use a fixed-size sieve region $H = [-B,B[ \times [-B,B[ \times [0,B[$ 
where each lattice point will correspond to a norm
which is always divisible by $q$ and hopefully divisible by many $p$.  
Define lattices $\Lambda_q$ and $\Lambda_{pq}$ by
\begin{align*} L_q =
\begin{bmatrix}
q & -r & 0\\ 0 & 1 & -r\\ 0 & 0 & 1
\end{bmatrix}
,\hspace{5mm}
L_{pq} =
\begin{bmatrix}
pq & -t & 0\\ 0 & 1 & -t\\ 0 & 0 & 1
\end{bmatrix}
\end{align*}
where $f_0(r) \equiv 0 \mod{q}$ and $f_0(t) \equiv 0 \mod{pq}$,
and the columns are a basis. 
Compute an LLL-reduced basis for  both $\Lambda_q$ 
and $\Lambda_{pq}$ to get matrices $L_q'$ and $L_{pq}'$. Then let
\begin{align*}
L' = (L_q')^{-1} \cdot L_{pq}'
\end{align*}
which is an integer matrix by construction.
Let $\Lambda'$ be the lattice with basis $L'$.
We mark all $(i,j,k)$ in $H \cap \Lambda'$.
As a result, for
a sieve location $(i,j,k)$ that has been marked, if we let
\begin{align*}
\begin{bmatrix}
a \\ b \\ c
\end{bmatrix} = L_q' \cdot 
\begin{bmatrix}
i \\ j \\ k
\end{bmatrix} 
\end{align*}
then we know that the norm of $\langle a + b\theta + c\theta^2 \rangle$ is divisible by both $q$ and $p$.  

We compute
and reduce $L_q$ once per special-$\mathfrak{q}$, and compute $L_{pq}$ etc for each $p$.
We compute $(a,b,c)$ only for $(i,j,k)$ that have been marked for many $p$ 
(above a  pre-determined threshold).

\bigskip

Our new results have two aspects.  
First, in Section 3 we improve the speed of enumeration of points in dimensions
higher than two.  Second, in Section 4 we give a new way of avoiding cache locality issues by the use of
a histogram of lattice point hits.  This applies regardless of dimension.

\section{Faster enumeration}

In a lattice $\Lambda$ of rank $n$ recall that the $i$-th successive minimum is defined by
\[
\lambda_i(\Lambda)=\inf \{ r \in \mathbb{R} :   \dim (\text{span}(\Lambda \cap B_r))\ge i \}
\]
where $B_r=\{x \in \mathbb{R}^n : ||x||\le r\}$.
In particular, $\lambda_1(\Lambda)$ is the length of a shortest nonzero vector in $\Lambda$.
A basis $v_1, \ldots ,v_n$  for $\Lambda$ is said to be a Minkowski-reduced basis if,
for $k=1,2, \ldots ,n$, $v_k$ is the shortest lattice element that can be extended
to a basis with $v_1, \ldots ,v_{k-1}$.

We assume we are sieving in three dimensions.
We fix a bound $B$ and let
$$H = [-B,B[ \times [-B,B[ \times [0,B[$$ be the sieving region.
Let  $\Lambda'$
be the lattice defined in Section 2.
The problem is  to list the elements of $\Lambda' \cap H$  in an efficient way.
In previous work this is done by going through the planes parallel to the $xy$-plane, 
and enumerating the lattice points in each of these planes.
We propose a different method which uses fewer planes.

Let $v_1, v_2, v_3$ be vectors having lengths 
$\lambda_1 (\Lambda'), \lambda_2 (\Lambda'), \lambda_3 (\Lambda')$,
 the first three successive minima of $\Lambda'$.
These three vectors are guaranteed to exist 
and we can either find all three, or an acceptably close approximation (see Remark 
\ref{LLL}).
The origin together with $v_1$ and $v_2$ define a plane which we call $P$. 
Let 
\[
c_{max} = \text{max}\{c \in \N : H \cap \left( P - c \cdot v_3 \right) \neq \emptyset\}.
\]
We refer to the plane $G=P - c_{max} \cdot v_3$ as the `ground plane'.

Our approach is very simple: to enumerate all lattice points in $H$,
we enumerate all points in the ground plane $G$ that lie in $H$, 
and then all points in the translates $G+kv_3$ for 
$k=1,2,3,...$  that lie in $H$, until we reach the last translate intersecting $H$.

\begin{remark}\label{LLL}
Finding $v_1, v_2, v_3$ is done with the LLL algorithm.  In practice, in very small
dimension such as three, this is sufficient to find a Minkowski-reduced basis, or a close approximation
which is good enough for our purposes.
\end{remark}

\begin{remark}
To easily enumerate points in a plane $G+kv_3$, we first  locate one point 
$p_0$ that is contained within  the plane and 
the sieving region $H$.  For this, we use integer linear programming (described in this context below).
Once we have located $p_0$, we proceed to enumerate points in this plane by adding and
subtracting multiples of $v_1$ and $v_2$ from $p_0$, until by doing so we are no longer
within $H$.
This is done inductively, by  first enumerating all $p_0+c_1 v_1$ where $c_1$ runs over
all integers such that $p_0+c_1 v_1$ is in $H$.
Then we add $v_2$ and enumerate all $p_0+v_2 + c_1 v_1$ where $c_1$ runs over
all integers such that $p_0+v_2+c_1 v_1$ is in $H$. Then we add $v_2$ again, and repeat.
This may not be the optimal method of  enumerating points in $G+kv_3$, 
however it worked well in our computations
and is sufficient for our purposes.
Moreover, this inductive procedure will extend to higher dimensions,
as long as the integer linear programming problem required to find the corresponding
feasible points is tractable.  We expect this to be the case certainly up to dimension
six (which was previously thought to be out of reach) and perhaps further.
\end{remark}

\begin{remark}
If the lattice is very skewed, it is possible that the last valid sieving point
in the plane is $p_k = p_{k-1} + v_1 + c \cdot v_2$, where $c \ge 2$ 
(and $p_{k-1}$ is the previous point).
It would be preferable to be able to reach all points by unit additions of $v_2$
so for practical purposes, we do this and ignore the rare cases where such `outlier'
points are missed.
\end{remark}

\begin{remark}
In two dimensions, the sieving method of Franke and  Kleinjung \cite{SHARCS2005} is very efficient.
Our approach works in the 2d case also, using the first two successive minima of the 2d lattice, but it
will not quite compete with the method in \cite{SHARCS2005}
 in terms of speed of enumeration because we must do a little extra
work when dealing with boundaries.  This shows that dealing with the
boundaries of the sieving region  is not trivial.
\end{remark}

\subsection{Previous Lattice Enumeration Methods}

Lattice enumeration is widely used in algorithms to solve certain lattice problems, such
as the Closest Vector Problem.
However, sieving in a cuboid introduces many complications that do not occur when sieving in a ball.

Our enumeration here is similar to Babai's `nearest plane' algorithm for lattice enumeration
\cite{Babai}.  However, it is significantly different in that we sieve in a box,
as opposed to a sphere.  Further, we do not compute a norm for every point to test
if it is within the boundary - use of a box allows us to separate many points which
may be treated in fast loops with no individual boundary checking.  In practice
this makes a huge difference. Note that L. Gr\'emy's space sieve is 120 times
faster than Babai's algorithm (see \cite{GremyPhD}).
We outperform the space sieve by over $2.5\times$. Note also that sieving in a 
rectangular region is fundamental to Franke and Kleinjung's 2d lattice sieve algorithm,
and its success depends on the shape of this region.

\subsection{Lattice Width}

Our idea is to cover the lattice with as few hyperplanes as possible. 
This is motivated by the concept of `lattice width' which we now define.

Suppose we have a finite set of points $K\subseteq \mathbb{Z}^n$. 
Pick some direction $c\in \R^n$.
The width of $K$ in direction $c$ is defined to be
\[
w(K,c):=\sup_{x\in K} \langle c,x\rangle - \inf_{x\in K} \langle c,x\rangle
\]
where we only consider $c$ such that the supremum and infimum are finite.
Since $K$ is finite we can replace $\sup$ by $\max$ and $\inf$ by $\min$.
Geometrically, if $K$ has width $\ell$ in the direction  $c$ then
any element of $K$ lies on a hyperplane $ \langle c,x\rangle=b$
where $b$ is an integer between $\inf_{x\in K} \langle c,x\rangle$
and $\sup_{x\in K} \langle c,x\rangle$.
The idea is that $K$ has width $\ell$ in the direction $c$ if $K$ can be covered by 
$\lfloor \ell \rfloor +1$ parallel hyperplanes  which are orthogonal to $c$.

The lattice width of $K$ is defined to be the infimum
of the widths in the direction $c$, over all nonzero $c$ in the integer lattice:
\[
w(K):=\inf_{c\in \mathbb{Z}^n, c\not=0} w(K,c).
\]
Note that the lattice width is an integer because  $K\subseteq \mathbb{Z}^n$. 
Therefore, the lattice width tells us the minimal number of 
Diophantine hyperplanes needed to cover $K$.

There are techniques for calculating the lattice width, and the directions
that give it, however these are generally 
used for lattices in a high number of dimensions. Because we are only in three
dimensions our method of
using shortest vectors is simpler and is sufficient for our purposes.

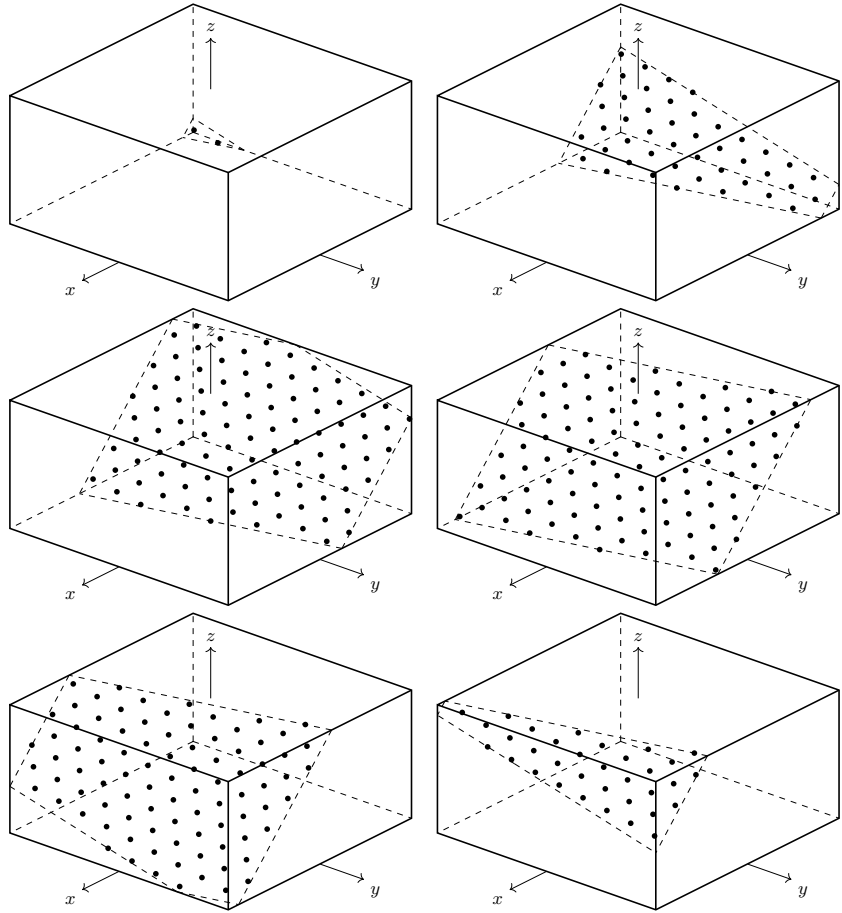
\begin{figure}[ht]
\centering
\scalebox{0.75}{%
\begin{tabular}{cc}
\tdplotsetmaincoords{65}{130}
\begin{tikzpicture}[tdplot_main_coords,scale=0.025]
\def\BigSide{100}
\def\SmallSide{60}
\pgfmathsetmacro{\CalcSide}{\BigSide-\SmallSide}

\tdplotsetcoord{P}{sqrt(3)*\BigSide}{55}{45}

\coordinate (sxl) at (\BigSide,\CalcSide,\BigSide);
\coordinate (syl) at (\CalcSide,\CalcSide,\BigSide);
\coordinate (szl) at (\CalcSide,\BigSide,\BigSide);

\draw[dashed] 
  (-100,-100,0) -- (-100,-100,100) 
  (-100,-100,0) -- (100,-100,0)
  (-100,-100,0) -- (-100,100,0);
\draw[->] 
  (Px) -- ++ (40,0,0) node[anchor=north east]{$x$};
\draw[->]
   (Py) -- ++(0,40,0) node[anchor=north west]{$y$};
\draw[->] 
  (Pz) -- ++(0,0,40) node[anchor=south]{$z$};

\draw[thick]
  (100,100,100) -- (100,-100,100) -- (100,-100,0) -- (100,100,0) -- (-100,100,0) --
  (-100,100,100) -- (-100,-100,100) -- (100,-100,100);
\draw[thick]
  (-100,100,100) -- (100,100,100) -- (100,100,0);

\node[fill,circle,inner sep=1.00pt] at (-96,-74,1) {};
\node[fill,circle,inner sep=1.00pt] at (-95,-95,5) {};
\draw[dashed]
   (-100,-100,100/9) -- (-3460/39,-100,0);
\draw[dashed]
   (-100,-2500/47,0) -- (-3460/39,-100,0);
\draw[dashed]
   (-100,-2500/47,0) -- (-100,-100,100/9);

\end{tikzpicture} &
\tdplotsetmaincoords{65}{130}
\begin{tikzpicture}[tdplot_main_coords,scale=0.025]
\def\BigSide{100}
\def\SmallSide{60}
\pgfmathsetmacro{\CalcSide}{\BigSide-\SmallSide}

\tdplotsetcoord{P}{sqrt(3)*\BigSide}{55}{45}

\coordinate (sxl) at (\BigSide,\CalcSide,\BigSide);
\coordinate (syl) at (\CalcSide,\CalcSide,\BigSide);
\coordinate (szl) at (\CalcSide,\BigSide,\BigSide);

\draw[dashed] 
  (-100,-100,0) -- (-100,-100,100) 
  (-100,-100,0) -- (100,-100,0)
  (-100,-100,0) -- (-100,100,0);
\draw[->] 
  (Px) -- ++ (40,0,0) node[anchor=north east]{$x$};
\draw[->]
   (Py) -- ++(0,40,0) node[anchor=north west]{$y$};
\draw[->] 
  (Pz) -- ++(0,0,40) node[anchor=south]{$z$};

\draw[thick]
  (100,100,100) -- (100,-100,100) -- (100,-100,0) -- (100,100,0) -- (-100,100,0) --
  (-100,100,100) -- (-100,-100,100) -- (100,-100,100);
\draw[thick]
  (-100,100,100) -- (100,100,100) -- (100,100,0);

\node[fill,circle,inner sep=1.00pt] at (-87,89,9) {};
\node[fill,circle,inner sep=1.00pt] at (-77,77,2) {};
\node[fill,circle,inner sep=1.00pt] at (-96,80,20) {};
\node[fill,circle,inner sep=1.00pt] at (-86,68,13) {};
\node[fill,circle,inner sep=1.00pt] at (-76,56,6) {};
\node[fill,circle,inner sep=1.00pt] at (-95,59,24) {};
\node[fill,circle,inner sep=1.00pt] at (-85,47,17) {};
\node[fill,circle,inner sep=1.00pt] at (-75,35,10) {};
\node[fill,circle,inner sep=1.00pt] at (-65,23,3) {};
\node[fill,circle,inner sep=1.00pt] at (-94,38,28) {};
\node[fill,circle,inner sep=1.00pt] at (-84,26,21) {};
\node[fill,circle,inner sep=1.00pt] at (-74,14,14) {};
\node[fill,circle,inner sep=1.00pt] at (-64,2,7) {};
\node[fill,circle,inner sep=1.00pt] at (-54,-10,0) {};
\node[fill,circle,inner sep=1.00pt] at (-93,17,32) {};
\node[fill,circle,inner sep=1.00pt] at (-83,5,25) {};
\node[fill,circle,inner sep=1.00pt] at (-73,-7,18) {};
\node[fill,circle,inner sep=1.00pt] at (-63,-19,11) {};
\node[fill,circle,inner sep=1.00pt] at (-53,-31,4) {};
\node[fill,circle,inner sep=1.00pt] at (-92,-4,36) {};
\node[fill,circle,inner sep=1.00pt] at (-82,-16,29) {};
\node[fill,circle,inner sep=1.00pt] at (-72,-28,22) {};
\node[fill,circle,inner sep=1.00pt] at (-62,-40,15) {};
\node[fill,circle,inner sep=1.00pt] at (-52,-52,8) {};
\node[fill,circle,inner sep=1.00pt] at (-42,-64,1) {};
\node[fill,circle,inner sep=1.00pt] at (-91,-25,40) {};
\node[fill,circle,inner sep=1.00pt] at (-81,-37,33) {};
\node[fill,circle,inner sep=1.00pt] at (-71,-49,26) {};
\node[fill,circle,inner sep=1.00pt] at (-61,-61,19) {};
\node[fill,circle,inner sep=1.00pt] at (-51,-73,12) {};
\node[fill,circle,inner sep=1.00pt] at (-41,-85,5) {};
\node[fill,circle,inner sep=1.00pt] at (-100,-34,51) {};
\node[fill,circle,inner sep=1.00pt] at (-90,-46,44) {};
\node[fill,circle,inner sep=1.00pt] at (-80,-58,37) {};
\node[fill,circle,inner sep=1.00pt] at (-70,-70,30) {};
\node[fill,circle,inner sep=1.00pt] at (-60,-82,23) {};
\node[fill,circle,inner sep=1.00pt] at (-50,-94,16) {};
\node[fill,circle,inner sep=1.00pt] at (-99,-55,55) {};
\node[fill,circle,inner sep=1.00pt] at (-89,-67,48) {};
\node[fill,circle,inner sep=1.00pt] at (-79,-79,41) {};
\node[fill,circle,inner sep=1.00pt] at (-69,-91,34) {};
\node[fill,circle,inner sep=1.00pt] at (-98,-76,59) {};
\node[fill,circle,inner sep=1.00pt] at (-88,-88,52) {};
\node[fill,circle,inner sep=1.00pt] at (-78,-100,45) {};
\node[fill,circle,inner sep=1.00pt] at (-97,-97,63) {};
\draw[dashed]
   (-100,-100,200/3) -- (-420/13,-100,0);
\draw[dashed]
   (-3140/39,100,0) -- (-420/13,-100,0);
\draw[dashed]
   (-3140/39,100,0) -- (-100,100,1900/99);
\draw[dashed]
   (-100,-100,200/3) -- (-100,100,1900/99);
\end{tikzpicture} \\

\tdplotsetmaincoords{65}{130}
\begin{tikzpicture}[tdplot_main_coords,scale=0.025]
\def\BigSide{100}
\def\SmallSide{60}
\pgfmathsetmacro{\CalcSide}{\BigSide-\SmallSide}

\tdplotsetcoord{P}{sqrt(3)*\BigSide}{55}{45}

\coordinate (sxl) at (\BigSide,\CalcSide,\BigSide);
\coordinate (syl) at (\CalcSide,\CalcSide,\BigSide);
\coordinate (szl) at (\CalcSide,\BigSide,\BigSide);

\draw[dashed] 
  (-100,-100,0) -- (-100,-100,100) 
  (-100,-100,0) -- (100,-100,0)
  (-100,-100,0) -- (-100,100,0);
\draw[->] 
  (Px) -- ++ (40,0,0) node[anchor=north east]{$x$};
\draw[->]
   (Py) -- ++(0,40,0) node[anchor=north west]{$y$};
\draw[->] 
  (Pz) -- ++(0,0,40) node[anchor=south]{$z$};

\draw[thick]
  (100,100,100) -- (100,-100,100) -- (100,-100,0) -- (100,100,0) -- (-100,100,0) --
  (-100,100,100) -- (-100,-100,100) -- (100,-100,100);
\draw[thick]
  (-100,100,100) -- (100,100,100) -- (100,100,0);
\node[fill,circle,inner sep=1.00pt] at (-33,99,9) {};
\node[fill,circle,inner sep=1.00pt] at (-23,87,2) {};
\node[fill,circle,inner sep=1.00pt] at (-42,90,20) {};
\node[fill,circle,inner sep=1.00pt] at (-32,78,13) {};
\node[fill,circle,inner sep=1.00pt] at (-22,66,6) {};
\node[fill,circle,inner sep=1.00pt] at (-61,93,38) {};
\node[fill,circle,inner sep=1.00pt] at (-51,81,31) {};
\node[fill,circle,inner sep=1.00pt] at (-41,69,24) {};
\node[fill,circle,inner sep=1.00pt] at (-31,57,17) {};
\node[fill,circle,inner sep=1.00pt] at (-21,45,10) {};
\node[fill,circle,inner sep=1.00pt] at (-11,33,3) {};
\node[fill,circle,inner sep=1.00pt] at (-80,96,56) {};
\node[fill,circle,inner sep=1.00pt] at (-70,84,49) {};
\node[fill,circle,inner sep=1.00pt] at (-60,72,42) {};
\node[fill,circle,inner sep=1.00pt] at (-50,60,35) {};
\node[fill,circle,inner sep=1.00pt] at (-40,48,28) {};
\node[fill,circle,inner sep=1.00pt] at (-30,36,21) {};
\node[fill,circle,inner sep=1.00pt] at (-20,24,14) {};
\node[fill,circle,inner sep=1.00pt] at (-10,12,7) {};
\node[fill,circle,inner sep=1.00pt] at (0,0,0) {};
\node[fill,circle,inner sep=1.00pt] at (-99,99,74) {};
\node[fill,circle,inner sep=1.00pt] at (-89,87,67) {};
\node[fill,circle,inner sep=1.00pt] at (-79,75,60) {};
\node[fill,circle,inner sep=1.00pt] at (-69,63,53) {};
\node[fill,circle,inner sep=1.00pt] at (-59,51,46) {};
\node[fill,circle,inner sep=1.00pt] at (-49,39,39) {};
\node[fill,circle,inner sep=1.00pt] at (-39,27,32) {};
\node[fill,circle,inner sep=1.00pt] at (-29,15,25) {};
\node[fill,circle,inner sep=1.00pt] at (-19,3,18) {};
\node[fill,circle,inner sep=1.00pt] at (-9,-9,11) {};
\node[fill,circle,inner sep=1.00pt] at (1,-21,4) {};
\node[fill,circle,inner sep=1.00pt] at (-98,78,78) {};
\node[fill,circle,inner sep=1.00pt] at (-88,66,71) {};
\node[fill,circle,inner sep=1.00pt] at (-78,54,64) {};
\node[fill,circle,inner sep=1.00pt] at (-68,42,57) {};
\node[fill,circle,inner sep=1.00pt] at (-58,30,50) {};
\node[fill,circle,inner sep=1.00pt] at (-48,18,43) {};
\node[fill,circle,inner sep=1.00pt] at (-38,6,36) {};
\node[fill,circle,inner sep=1.00pt] at (-28,-6,29) {};
\node[fill,circle,inner sep=1.00pt] at (-18,-18,22) {};
\node[fill,circle,inner sep=1.00pt] at (-8,-30,15) {};
\node[fill,circle,inner sep=1.00pt] at (2,-42,8) {};
\node[fill,circle,inner sep=1.00pt] at (12,-54,1) {};
\node[fill,circle,inner sep=1.00pt] at (-97,57,82) {};
\node[fill,circle,inner sep=1.00pt] at (-87,45,75) {};
\node[fill,circle,inner sep=1.00pt] at (-77,33,68) {};
\node[fill,circle,inner sep=1.00pt] at (-67,21,61) {};
\node[fill,circle,inner sep=1.00pt] at (-57,9,54) {};
\node[fill,circle,inner sep=1.00pt] at (-47,-3,47) {};
\node[fill,circle,inner sep=1.00pt] at (-37,-15,40) {};
\node[fill,circle,inner sep=1.00pt] at (-27,-27,33) {};
\node[fill,circle,inner sep=1.00pt] at (-17,-39,26) {};
\node[fill,circle,inner sep=1.00pt] at (-7,-51,19) {};
\node[fill,circle,inner sep=1.00pt] at (3,-63,12) {};
\node[fill,circle,inner sep=1.00pt] at (13,-75,5) {};
\node[fill,circle,inner sep=1.00pt] at (-96,36,86) {};
\node[fill,circle,inner sep=1.00pt] at (-86,24,79) {};
\node[fill,circle,inner sep=1.00pt] at (-76,12,72) {};
\node[fill,circle,inner sep=1.00pt] at (-66,0,65) {};
\node[fill,circle,inner sep=1.00pt] at (-56,-12,58) {};
\node[fill,circle,inner sep=1.00pt] at (-46,-24,51) {};
\node[fill,circle,inner sep=1.00pt] at (-36,-36,44) {};
\node[fill,circle,inner sep=1.00pt] at (-26,-48,37) {};
\node[fill,circle,inner sep=1.00pt] at (-16,-60,30) {};
\node[fill,circle,inner sep=1.00pt] at (-6,-72,23) {};
\node[fill,circle,inner sep=1.00pt] at (4,-84,16) {};
\node[fill,circle,inner sep=1.00pt] at (14,-96,9) {};
\node[fill,circle,inner sep=1.00pt] at (-95,15,90) {};
\node[fill,circle,inner sep=1.00pt] at (-85,3,83) {};
\node[fill,circle,inner sep=1.00pt] at (-75,-9,76) {};
\node[fill,circle,inner sep=1.00pt] at (-65,-21,69) {};
\node[fill,circle,inner sep=1.00pt] at (-55,-33,62) {};
\node[fill,circle,inner sep=1.00pt] at (-45,-45,55) {};
\node[fill,circle,inner sep=1.00pt] at (-35,-57,48) {};
\node[fill,circle,inner sep=1.00pt] at (-25,-69,41) {};
\node[fill,circle,inner sep=1.00pt] at (-15,-81,34) {};
\node[fill,circle,inner sep=1.00pt] at (-5,-93,27) {};
\node[fill,circle,inner sep=1.00pt] at (-94,-6,94) {};
\node[fill,circle,inner sep=1.00pt] at (-84,-18,87) {};
\node[fill,circle,inner sep=1.00pt] at (-74,-30,80) {};
\node[fill,circle,inner sep=1.00pt] at (-64,-42,73) {};
\node[fill,circle,inner sep=1.00pt] at (-54,-54,66) {};
\node[fill,circle,inner sep=1.00pt] at (-44,-66,59) {};
\node[fill,circle,inner sep=1.00pt] at (-34,-78,52) {};
\node[fill,circle,inner sep=1.00pt] at (-24,-90,45) {};
\node[fill,circle,inner sep=1.00pt] at (-93,-27,98) {};
\node[fill,circle,inner sep=1.00pt] at (-83,-39,91) {};
\node[fill,circle,inner sep=1.00pt] at (-73,-51,84) {};
\node[fill,circle,inner sep=1.00pt] at (-63,-63,77) {};
\node[fill,circle,inner sep=1.00pt] at (-53,-75,70) {};
\node[fill,circle,inner sep=1.00pt] at (-43,-87,63) {};
\node[fill,circle,inner sep=1.00pt] at (-33,-99,56) {};
\node[fill,circle,inner sep=1.00pt] at (-82,-60,95) {};
\node[fill,circle,inner sep=1.00pt] at (-72,-72,88) {};
\node[fill,circle,inner sep=1.00pt] at (-62,-84,81) {};
\node[fill,circle,inner sep=1.00pt] at (-52,-96,74) {};
\node[fill,circle,inner sep=1.00pt] at (-81,-81,99) {};
\node[fill,circle,inner sep=1.00pt] at (-71,-93,92) {};
\draw[dashed]
   (-3020/39,-100,100) -- (940/39,-100,0);
\draw[dashed]
   (-940/39,100,0) -- (940/39,-100,0);
\draw[dashed]
   (-940/39,100,0) -- (-100,100,7400/99);
\draw[dashed]
   (-100,-300/47,100) -- (-100,100,7400/99);
\draw[dashed]
   (-100,-300/47,100) -- (-3020/39,-100,100);
\end{tikzpicture} &
\tdplotsetmaincoords{65}{130}
\begin{tikzpicture}[tdplot_main_coords,scale=0.025]
\def\BigSide{100}
\def\SmallSide{60}
\pgfmathsetmacro{\CalcSide}{\BigSide-\SmallSide}

\tdplotsetcoord{P}{sqrt(3)*\BigSide}{55}{45}

\coordinate (sxl) at (\BigSide,\CalcSide,\BigSide);
\coordinate (syl) at (\CalcSide,\CalcSide,\BigSide);
\coordinate (szl) at (\CalcSide,\BigSide,\BigSide);

\draw[dashed] 
  (-100,-100,0) -- (-100,-100,100) 
  (-100,-100,0) -- (100,-100,0)
  (-100,-100,0) -- (-100,100,0);
\draw[->] 
  (Px) -- ++ (40,0,0) node[anchor=north east]{$x$};
\draw[->]
   (Py) -- ++(0,40,0) node[anchor=north west]{$y$};
\draw[->] 
  (Pz) -- ++(0,0,40) node[anchor=south]{$z$};

\draw[thick]
  (100,100,100) -- (100,-100,100) -- (100,-100,0) -- (100,100,0) -- (-100,100,0) --
  (-100,100,100) -- (-100,-100,100) -- (100,-100,100);
\draw[thick]
  (-100,100,100) -- (100,100,100) -- (100,100,0);

\node[fill,circle,inner sep=1.00pt] at (31,97,2) {};
\node[fill,circle,inner sep=1.00pt] at (12,100,20) {};
\node[fill,circle,inner sep=1.00pt] at (22,88,13) {};
\node[fill,circle,inner sep=1.00pt] at (32,76,6) {};
\node[fill,circle,inner sep=1.00pt] at (3,91,31) {};
\node[fill,circle,inner sep=1.00pt] at (13,79,24) {};
\node[fill,circle,inner sep=1.00pt] at (23,67,17) {};
\node[fill,circle,inner sep=1.00pt] at (33,55,10) {};
\node[fill,circle,inner sep=1.00pt] at (43,43,3) {};
\node[fill,circle,inner sep=1.00pt] at (-16,94,49) {};
\node[fill,circle,inner sep=1.00pt] at (-6,82,42) {};
\node[fill,circle,inner sep=1.00pt] at (4,70,35) {};
\node[fill,circle,inner sep=1.00pt] at (14,58,28) {};
\node[fill,circle,inner sep=1.00pt] at (24,46,21) {};
\node[fill,circle,inner sep=1.00pt] at (34,34,14) {};
\node[fill,circle,inner sep=1.00pt] at (44,22,7) {};
\node[fill,circle,inner sep=1.00pt] at (54,10,0) {};
\node[fill,circle,inner sep=1.00pt] at (-35,97,67) {};
\node[fill,circle,inner sep=1.00pt] at (-25,85,60) {};
\node[fill,circle,inner sep=1.00pt] at (-15,73,53) {};
\node[fill,circle,inner sep=1.00pt] at (-5,61,46) {};
\node[fill,circle,inner sep=1.00pt] at (5,49,39) {};
\node[fill,circle,inner sep=1.00pt] at (15,37,32) {};
\node[fill,circle,inner sep=1.00pt] at (25,25,25) {};
\node[fill,circle,inner sep=1.00pt] at (35,13,18) {};
\node[fill,circle,inner sep=1.00pt] at (45,1,11) {};
\node[fill,circle,inner sep=1.00pt] at (55,-11,4) {};
\node[fill,circle,inner sep=1.00pt] at (-54,100,85) {};
\node[fill,circle,inner sep=1.00pt] at (-44,88,78) {};
\node[fill,circle,inner sep=1.00pt] at (-34,76,71) {};
\node[fill,circle,inner sep=1.00pt] at (-24,64,64) {};
\node[fill,circle,inner sep=1.00pt] at (-14,52,57) {};
\node[fill,circle,inner sep=1.00pt] at (-4,40,50) {};
\node[fill,circle,inner sep=1.00pt] at (6,28,43) {};
\node[fill,circle,inner sep=1.00pt] at (16,16,36) {};
\node[fill,circle,inner sep=1.00pt] at (26,4,29) {};
\node[fill,circle,inner sep=1.00pt] at (36,-8,22) {};
\node[fill,circle,inner sep=1.00pt] at (46,-20,15) {};
\node[fill,circle,inner sep=1.00pt] at (56,-32,8) {};
\node[fill,circle,inner sep=1.00pt] at (66,-44,1) {};
\node[fill,circle,inner sep=1.00pt] at (-63,91,96) {};
\node[fill,circle,inner sep=1.00pt] at (-53,79,89) {};
\node[fill,circle,inner sep=1.00pt] at (-43,67,82) {};
\node[fill,circle,inner sep=1.00pt] at (-33,55,75) {};
\node[fill,circle,inner sep=1.00pt] at (-23,43,68) {};
\node[fill,circle,inner sep=1.00pt] at (-13,31,61) {};
\node[fill,circle,inner sep=1.00pt] at (-3,19,54) {};
\node[fill,circle,inner sep=1.00pt] at (7,7,47) {};
\node[fill,circle,inner sep=1.00pt] at (17,-5,40) {};
\node[fill,circle,inner sep=1.00pt] at (27,-17,33) {};
\node[fill,circle,inner sep=1.00pt] at (37,-29,26) {};
\node[fill,circle,inner sep=1.00pt] at (47,-41,19) {};
\node[fill,circle,inner sep=1.00pt] at (57,-53,12) {};
\node[fill,circle,inner sep=1.00pt] at (67,-65,5) {};
\node[fill,circle,inner sep=1.00pt] at (-62,70,100) {};
\node[fill,circle,inner sep=1.00pt] at (-52,58,93) {};
\node[fill,circle,inner sep=1.00pt] at (-42,46,86) {};
\node[fill,circle,inner sep=1.00pt] at (-32,34,79) {};
\node[fill,circle,inner sep=1.00pt] at (-22,22,72) {};
\node[fill,circle,inner sep=1.00pt] at (-12,10,65) {};
\node[fill,circle,inner sep=1.00pt] at (-2,-2,58) {};
\node[fill,circle,inner sep=1.00pt] at (8,-14,51) {};
\node[fill,circle,inner sep=1.00pt] at (18,-26,44) {};
\node[fill,circle,inner sep=1.00pt] at (28,-38,37) {};
\node[fill,circle,inner sep=1.00pt] at (38,-50,30) {};
\node[fill,circle,inner sep=1.00pt] at (48,-62,23) {};
\node[fill,circle,inner sep=1.00pt] at (58,-74,16) {};
\node[fill,circle,inner sep=1.00pt] at (68,-86,9) {};
\node[fill,circle,inner sep=1.00pt] at (78,-98,2) {};
\node[fill,circle,inner sep=1.00pt] at (-51,37,97) {};
\node[fill,circle,inner sep=1.00pt] at (-41,25,90) {};
\node[fill,circle,inner sep=1.00pt] at (-31,13,83) {};
\node[fill,circle,inner sep=1.00pt] at (-21,1,76) {};
\node[fill,circle,inner sep=1.00pt] at (-11,-11,69) {};
\node[fill,circle,inner sep=1.00pt] at (-1,-23,62) {};
\node[fill,circle,inner sep=1.00pt] at (9,-35,55) {};
\node[fill,circle,inner sep=1.00pt] at (19,-47,48) {};
\node[fill,circle,inner sep=1.00pt] at (29,-59,41) {};
\node[fill,circle,inner sep=1.00pt] at (39,-71,34) {};
\node[fill,circle,inner sep=1.00pt] at (49,-83,27) {};
\node[fill,circle,inner sep=1.00pt] at (59,-95,20) {};
\node[fill,circle,inner sep=1.00pt] at (-40,4,94) {};
\node[fill,circle,inner sep=1.00pt] at (-30,-8,87) {};
\node[fill,circle,inner sep=1.00pt] at (-20,-20,80) {};
\node[fill,circle,inner sep=1.00pt] at (-10,-32,73) {};
\node[fill,circle,inner sep=1.00pt] at (0,-44,66) {};
\node[fill,circle,inner sep=1.00pt] at (10,-56,59) {};
\node[fill,circle,inner sep=1.00pt] at (20,-68,52) {};
\node[fill,circle,inner sep=1.00pt] at (30,-80,45) {};
\node[fill,circle,inner sep=1.00pt] at (40,-92,38) {};
\node[fill,circle,inner sep=1.00pt] at (-39,-17,98) {};
\node[fill,circle,inner sep=1.00pt] at (-29,-29,91) {};
\node[fill,circle,inner sep=1.00pt] at (-19,-41,84) {};
\node[fill,circle,inner sep=1.00pt] at (-9,-53,77) {};
\node[fill,circle,inner sep=1.00pt] at (1,-65,70) {};
\node[fill,circle,inner sep=1.00pt] at (11,-77,63) {};
\node[fill,circle,inner sep=1.00pt] at (21,-89,56) {};
\node[fill,circle,inner sep=1.00pt] at (-28,-50,95) {};
\node[fill,circle,inner sep=1.00pt] at (-18,-62,88) {};
\node[fill,circle,inner sep=1.00pt] at (-8,-74,81) {};
\node[fill,circle,inner sep=1.00pt] at (2,-86,74) {};
\node[fill,circle,inner sep=1.00pt] at (12,-98,67) {};
\node[fill,circle,inner sep=1.00pt] at (-27,-71,99) {};
\node[fill,circle,inner sep=1.00pt] at (-17,-83,92) {};
\node[fill,circle,inner sep=1.00pt] at (-7,-95,85) {};
\draw[dashed]
   (-820/39,-100,100) -- (3140/39,-100,0);
\draw[dashed]
   (420/13,100,0) -- (3140/39,-100,0);
\draw[dashed]
   (420/13,100,0) -- (-900/13,100,100);
\draw[dashed]
   (-820/39,-100,100) -- (-900/13,100,100);

\end{tikzpicture}\\

\tdplotsetmaincoords{65}{130}
\begin{tikzpicture}[tdplot_main_coords,scale=0.025]
\def\BigSide{100}
\def\SmallSide{60}
\pgfmathsetmacro{\CalcSide}{\BigSide-\SmallSide}

\tdplotsetcoord{P}{sqrt(3)*\BigSide}{55}{45}

\coordinate (sxl) at (\BigSide,\CalcSide,\BigSide);
\coordinate (syl) at (\CalcSide,\CalcSide,\BigSide);
\coordinate (szl) at (\CalcSide,\BigSide,\BigSide);

\draw[dashed] 
  (-100,-100,0) -- (-100,-100,100) 
  (-100,-100,0) -- (100,-100,0)
  (-100,-100,0) -- (-100,100,0);
\draw[->] 
  (Px) -- ++ (40,0,0) node[anchor=north east]{$x$};
\draw[->]
   (Py) -- ++(0,40,0) node[anchor=north west]{$y$};
\draw[->] 
  (Pz) -- ++(0,0,40) node[anchor=south]{$z$};

\draw[thick]
  (100,100,100) -- (100,-100,100) -- (100,-100,0) -- (100,100,0) -- (-100,100,0) --
  (-100,100,100) -- (-100,-100,100) -- (100,-100,100);
\draw[thick]
  (-100,100,100) -- (100,100,100) -- (100,100,0);

\node[fill,circle,inner sep=1.00pt] at (76,98,13) {};
\node[fill,circle,inner sep=1.00pt] at (86,86,6) {};
\node[fill,circle,inner sep=1.00pt] at (67,89,24) {};
\node[fill,circle,inner sep=1.00pt] at (77,77,17) {};
\node[fill,circle,inner sep=1.00pt] at (87,65,10) {};
\node[fill,circle,inner sep=1.00pt] at (97,53,3) {};
\node[fill,circle,inner sep=1.00pt] at (48,92,42) {};
\node[fill,circle,inner sep=1.00pt] at (58,80,35) {};
\node[fill,circle,inner sep=1.00pt] at (68,68,28) {};
\node[fill,circle,inner sep=1.00pt] at (78,56,21) {};
\node[fill,circle,inner sep=1.00pt] at (88,44,14) {};
\node[fill,circle,inner sep=1.00pt] at (98,32,7) {};
\node[fill,circle,inner sep=1.00pt] at (29,95,60) {};
\node[fill,circle,inner sep=1.00pt] at (39,83,53) {};
\node[fill,circle,inner sep=1.00pt] at (49,71,46) {};
\node[fill,circle,inner sep=1.00pt] at (59,59,39) {};
\node[fill,circle,inner sep=1.00pt] at (69,47,32) {};
\node[fill,circle,inner sep=1.00pt] at (79,35,25) {};
\node[fill,circle,inner sep=1.00pt] at (89,23,18) {};
\node[fill,circle,inner sep=1.00pt] at (99,11,11) {};
\node[fill,circle,inner sep=1.00pt] at (10,98,78) {};
\node[fill,circle,inner sep=1.00pt] at (20,86,71) {};
\node[fill,circle,inner sep=1.00pt] at (30,74,64) {};
\node[fill,circle,inner sep=1.00pt] at (40,62,57) {};
\node[fill,circle,inner sep=1.00pt] at (50,50,50) {};
\node[fill,circle,inner sep=1.00pt] at (60,38,43) {};
\node[fill,circle,inner sep=1.00pt] at (70,26,36) {};
\node[fill,circle,inner sep=1.00pt] at (80,14,29) {};
\node[fill,circle,inner sep=1.00pt] at (90,2,22) {};
\node[fill,circle,inner sep=1.00pt] at (100,-10,15) {};
\node[fill,circle,inner sep=1.00pt] at (1,89,89) {};
\node[fill,circle,inner sep=1.00pt] at (11,77,82) {};
\node[fill,circle,inner sep=1.00pt] at (21,65,75) {};
\node[fill,circle,inner sep=1.00pt] at (31,53,68) {};
\node[fill,circle,inner sep=1.00pt] at (41,41,61) {};
\node[fill,circle,inner sep=1.00pt] at (51,29,54) {};
\node[fill,circle,inner sep=1.00pt] at (61,17,47) {};
\node[fill,circle,inner sep=1.00pt] at (71,5,40) {};
\node[fill,circle,inner sep=1.00pt] at (81,-7,33) {};
\node[fill,circle,inner sep=1.00pt] at (91,-19,26) {};
\node[fill,circle,inner sep=1.00pt] at (-8,80,100) {};
\node[fill,circle,inner sep=1.00pt] at (2,68,93) {};
\node[fill,circle,inner sep=1.00pt] at (12,56,86) {};
\node[fill,circle,inner sep=1.00pt] at (22,44,79) {};
\node[fill,circle,inner sep=1.00pt] at (32,32,72) {};
\node[fill,circle,inner sep=1.00pt] at (42,20,65) {};
\node[fill,circle,inner sep=1.00pt] at (52,8,58) {};
\node[fill,circle,inner sep=1.00pt] at (62,-4,51) {};
\node[fill,circle,inner sep=1.00pt] at (72,-16,44) {};
\node[fill,circle,inner sep=1.00pt] at (82,-28,37) {};
\node[fill,circle,inner sep=1.00pt] at (92,-40,30) {};
\node[fill,circle,inner sep=1.00pt] at (3,47,97) {};
\node[fill,circle,inner sep=1.00pt] at (13,35,90) {};
\node[fill,circle,inner sep=1.00pt] at (23,23,83) {};
\node[fill,circle,inner sep=1.00pt] at (33,11,76) {};
\node[fill,circle,inner sep=1.00pt] at (43,-1,69) {};
\node[fill,circle,inner sep=1.00pt] at (53,-13,62) {};
\node[fill,circle,inner sep=1.00pt] at (63,-25,55) {};
\node[fill,circle,inner sep=1.00pt] at (73,-37,48) {};
\node[fill,circle,inner sep=1.00pt] at (83,-49,41) {};
\node[fill,circle,inner sep=1.00pt] at (93,-61,34) {};
\node[fill,circle,inner sep=1.00pt] at (14,14,94) {};
\node[fill,circle,inner sep=1.00pt] at (24,2,87) {};
\node[fill,circle,inner sep=1.00pt] at (34,-10,80) {};
\node[fill,circle,inner sep=1.00pt] at (44,-22,73) {};
\node[fill,circle,inner sep=1.00pt] at (54,-34,66) {};
\node[fill,circle,inner sep=1.00pt] at (64,-46,59) {};
\node[fill,circle,inner sep=1.00pt] at (74,-58,52) {};
\node[fill,circle,inner sep=1.00pt] at (84,-70,45) {};
\node[fill,circle,inner sep=1.00pt] at (94,-82,38) {};
\node[fill,circle,inner sep=1.00pt] at (15,-7,98) {};
\node[fill,circle,inner sep=1.00pt] at (25,-19,91) {};
\node[fill,circle,inner sep=1.00pt] at (35,-31,84) {};
\node[fill,circle,inner sep=1.00pt] at (45,-43,77) {};
\node[fill,circle,inner sep=1.00pt] at (55,-55,70) {};
\node[fill,circle,inner sep=1.00pt] at (65,-67,63) {};
\node[fill,circle,inner sep=1.00pt] at (75,-79,56) {};
\node[fill,circle,inner sep=1.00pt] at (85,-91,49) {};
\node[fill,circle,inner sep=1.00pt] at (26,-40,95) {};
\node[fill,circle,inner sep=1.00pt] at (36,-52,88) {};
\node[fill,circle,inner sep=1.00pt] at (46,-64,81) {};
\node[fill,circle,inner sep=1.00pt] at (56,-76,74) {};
\node[fill,circle,inner sep=1.00pt] at (66,-88,67) {};
\node[fill,circle,inner sep=1.00pt] at (76,-100,60) {};
\node[fill,circle,inner sep=1.00pt] at (27,-61,99) {};
\node[fill,circle,inner sep=1.00pt] at (37,-73,92) {};
\node[fill,circle,inner sep=1.00pt] at (47,-85,85) {};
\node[fill,circle,inner sep=1.00pt] at (57,-97,78) {};
\node[fill,circle,inner sep=1.00pt] at (38,-94,96) {};
\draw[dashed]
   (460/13,-100,100) -- (100,-100,400/11);
\draw[dashed]
   (100,2500/47,0) -- (100,-100,400/11);
\draw[dashed]
   (100,2500/47,0) -- (3460/39,100,0);
\draw[dashed]
   (-500/39,100,100) -- (3460/39,100,0);
\draw[dashed]
   (-500/39,100,100) -- (460/13,-100,100);

\end{tikzpicture} &
\tdplotsetmaincoords{65}{130}
\begin{tikzpicture}[tdplot_main_coords,scale=0.025]
\def\BigSide{100}
\def\SmallSide{60}
\pgfmathsetmacro{\CalcSide}{\BigSide-\SmallSide}

\tdplotsetcoord{P}{sqrt(3)*\BigSide}{55}{45}

\coordinate (sxl) at (\BigSide,\CalcSide,\BigSide);
\coordinate (syl) at (\CalcSide,\CalcSide,\BigSide);
\coordinate (szl) at (\CalcSide,\BigSide,\BigSide);

\draw[dashed] 
  (-100,-100,0) -- (-100,-100,100) 
  (-100,-100,0) -- (100,-100,0)
  (-100,-100,0) -- (-100,100,0);
\draw[->] 
  (Px) -- ++ (40,0,0) node[anchor=north east]{$x$};
\draw[->]
   (Py) -- ++(0,40,0) node[anchor=north west]{$y$};
\draw[->] 
  (Pz) -- ++(0,0,40) node[anchor=south]{$z$};

\draw[thick]
  (100,100,100) -- (100,-100,100) -- (100,-100,0) -- (100,100,0) -- (-100,100,0) --
  (-100,100,100) -- (-100,-100,100) -- (100,-100,100);
\draw[thick]
  (-100,100,100) -- (100,100,100) -- (100,100,0);

\node[fill,circle,inner sep=1.00pt] at (93,93,53) {};
\node[fill,circle,inner sep=1.00pt] at (74,96,71) {};
\node[fill,circle,inner sep=1.00pt] at (84,84,64) {};
\node[fill,circle,inner sep=1.00pt] at (94,72,57) {};
\node[fill,circle,inner sep=1.00pt] at (55,99,89) {};
\node[fill,circle,inner sep=1.00pt] at (65,87,82) {};
\node[fill,circle,inner sep=1.00pt] at (75,75,75) {};
\node[fill,circle,inner sep=1.00pt] at (85,63,68) {};
\node[fill,circle,inner sep=1.00pt] at (95,51,61) {};
\node[fill,circle,inner sep=1.00pt] at (46,90,100) {};
\node[fill,circle,inner sep=1.00pt] at (56,78,93) {};
\node[fill,circle,inner sep=1.00pt] at (66,66,86) {};
\node[fill,circle,inner sep=1.00pt] at (76,54,79) {};
\node[fill,circle,inner sep=1.00pt] at (86,42,72) {};
\node[fill,circle,inner sep=1.00pt] at (96,30,65) {};
\node[fill,circle,inner sep=1.00pt] at (57,57,97) {};
\node[fill,circle,inner sep=1.00pt] at (67,45,90) {};
\node[fill,circle,inner sep=1.00pt] at (77,33,83) {};
\node[fill,circle,inner sep=1.00pt] at (87,21,76) {};
\node[fill,circle,inner sep=1.00pt] at (97,9,69) {};
\node[fill,circle,inner sep=1.00pt] at (68,24,94) {};
\node[fill,circle,inner sep=1.00pt] at (78,12,87) {};
\node[fill,circle,inner sep=1.00pt] at (88,0,80) {};
\node[fill,circle,inner sep=1.00pt] at (98,-12,73) {};
\node[fill,circle,inner sep=1.00pt] at (69,3,98) {};
\node[fill,circle,inner sep=1.00pt] at (79,-9,91) {};
\node[fill,circle,inner sep=1.00pt] at (89,-21,84) {};
\node[fill,circle,inner sep=1.00pt] at (99,-33,77) {};
\node[fill,circle,inner sep=1.00pt] at (80,-30,95) {};
\node[fill,circle,inner sep=1.00pt] at (90,-42,88) {};
\node[fill,circle,inner sep=1.00pt] at (100,-54,81) {};
\node[fill,circle,inner sep=1.00pt] at (81,-51,99) {};
\node[fill,circle,inner sep=1.00pt] at (91,-63,92) {};
\node[fill,circle,inner sep=1.00pt] at (92,-84,96) {};
\draw[dashed]
   (3580/39,-100,100) -- (100,-100,9100/99);
\draw[dashed]
   (100,-500/17,11500/153) -- (100,-100,9100/99);
\draw[dashed]
   (100,-500/17,11500/153) -- (100,100,400/9);
\draw[dashed]
   (1700/39,100,100) -- (100,100,400/9);
\draw[dashed]
   (1700/39,100,100) -- (3580/39,-100,100);

\end{tikzpicture}
\end{tabular}
}
\caption{Six dense sublattices cover every point in the sieve region.}
\label{fig:sympathetic}
\end{figure}

\bigskip
\noindent {\bf Example} 
The lattice used in Fig. \ref{fig:sympathetic} is the following:
\begin{align*}
\begin{bmatrix}
10 & 18 & 35\\ -12 & 18 & 13\\ -7 & -22 & 18
\end{bmatrix}
\end{align*}
 The sieve region is $[-100,100]\times[-100,100]\times[0,100]$.  
 In this example, using our method,  6 planes cover every valid
sieving point.  With traditional plane sieving, using planes that are parallel to the base of the sieving cuboid, 101 planes are needed, each with at most  four lattice points.

\subsection{Integer Linear Programming}
Given a plane defined by  $v_1,v_2$ and a point $R$, with $R$ not necessarily
contained in the sieving region defined by $H = [-B,B[\times[-B,B[\times[0,B[$, the task is to find a point $p_0 = (x_0,y_0,z_0)$ that
is provably contained in the intersection of the plane and $H$, if such
a point exists.  We look for $r,s \in \Z$ such that $p_0 = R + r\cdot v_1 + s\cdot v_2$
and $p_0\in H$.

\noindent This can be formulated as an integer linear programming problem, where the
aim is to minimize $x$, subject to
\begin{flalign*}
A \cdot x \leq b
\end{flalign*}
where $x = (r,s) \in \Z^2$ and $A \in M_2(\Z), b \in \Z^2$, depend on $v_1,v_2,B$, and we
must find any feasible point, if one exists.  This problem is well studied, and though
it is NP-hard in general, can be solved easily in small dimensions.  It is
computationally trivial in dimension 3, for example, which we use in this article.

\section{Improved Cache Locality}\label{histogram}

Representing a lattice in memory is not necessarily done best using the `obvious' approach
of arranging all possible element co-ordinates in lines/planes and so on, and
then accessing points via a canonical list of co-ordinate places.  The
storage/retrieval of points tends to result in random memory access patterns, which
severely impacts performance.  This is a fundamental concern in large-scale computation.
Computer manufacturers address this by providing various levels of `cache', i.e. a limited
quantity of high-speed memory, too costly for main memory, which is used as a temporary
store of frequently-accessed or burst-access data.  The prior art in lattice-sieving
has always had to make use of cache to minimize the cost of the random memory
access patterns that occur in practice.

Our idea to improve cache locality  is simple: list and sort. 
We propose  to store all enumerated point coordinates in a list of increasing size.
Because the list increases strictly linearly, this is ideally suited to fast memory
access and is compatible with all levels of cache.  By itself, this is not an advantage
as we have merely collected a long list of randomly-organized points.  However, if we
encode points as e.g. a 32-bit integer, we can sort this list using these integers as a 
key.  Then, repeatedly-marked lattice points correspond to runs of identical keys.
If (key,value) pairs consist of such a key and a byte representing $\log{p}$, we
can recover lattice vectors with a large smooth part via a linear scan of the sorted
list.

Sorting is fast. When we consider that nowadays it is 
possible to sort one billion (key,value) integer pairs in seconds on a modern CPU,
it is evident that sorting, with its $O(N\log{N})$ or better complexity, is quite compatible
with modern cache hierarchies.  The situation is probably even better on GPU, although
it should be emphasized that in modern clusters, on GPU nodes there are typically
one or two GPUs and dozens of traditional CPUs, so it is not a priori obvious that one or
the other is to be preferred.

We compare sieving statistics in table \ref{SieveCompare} between our
implementation and that of \cite{MR3775580}.  Note that all of these times
give total special-$\mathfrak{q}$ time excluding the cofactorization time.
We have included listing/sorting times in our case.  In \cite{MR3775580},
sieving and memory access are intertwined and we compare this to our combined
sieving/listing/sorting time.  Cofactorization times are similar between the two.

\begin{table}[h]
\begin{center}
\begin{tabulary}{12cm}{|C|C|C|C|C|C|C|C|C|}
\hline
Authors&fbb&H&$q_{min}$&$q_{max}$&$\#\{q\}$&av. time&min time&max time\\
\hline
GGMT&$2^{21}$&10,10,8&16000000&16001000&7&143.93&142.17&145.28\\
GGMT&$2^{21}$&10,10,8&86500000&86501000&14&142.07&140.53&143.82\\
GGMT&$2^{21}$&10,10,8&262000000&262001000&9&142.40&140.95&144.34\\
\hline
GGMT&$2^{22}$&10,10,8&16000000&16001000&7&169.74&166.12&171.82\\
GGMT&$2^{22}$&10,10,8&86500000&86501000&14&167.53&166.01&173.55\\
GGMT&$2^{22}$&10,10,8&262000000&262001000&9&167.50&165.02&172.17\\
\hline
this work&$2^{24}$&9,10,10&16000000&16001000&7&35.47&34.94&36.95\\
this work&$2^{24}$&9,10,10&86500000&86501000&14&35.80&35.39&37.16\\
this work&$2^{24}$&9,10,10&262000000&262001000&9&36.37&35.71&37.54\\
\hline
\end{tabulary}
\end{center}
\caption{Sieve performance comparison (times in seconds)}\label{SieveCompare}
\end{table}

\section{Record computation in $F_{{p}^6}$}
\noindent We implemented the 3d case of our lattice sieving idea in C and used it to set a new record
in solving the discrete log in the multiplicative subgroup $(\F_{p^6})^\times$.  Previous records
were set by Zajac \cite{ZajacPhD}, Hayasaka et al (HAKT) \cite{HAKT}, and 
Gr\'{e}my et al (GGMT) \cite{MR3775580}.  All computations were done on the main compute
nodes of the Kay cluster at ICHEC, the Irish Center for High-End Computing.  Each node
consists of $2\times 20$ Intel Xeon Gold 6148 (Skylake) processors @ 2.4 GHz, with 192Gb RAM
per node.   All timings have been normalized to a nominal 2.0GHz clock speed.

With $\phi = (1 + \sqrt{5})/2$, we chose the prime $p = \floor{10^{21}\cdot \phi} + 29$.  Our target field
is $\F_{p^6}$, where $p^6$ has 423 bits.  This is comparable to the field size of the previous record at
422 bits \cite{MR3775580}.  One consequence of our choice is to allow a fair comparison of the total effort required to solve
discrete logs in a field of this order of magnitude.

\subsection{Polynomial selection}
We implemented the Joux-Lercier-Smart-Vercauteren ($\text{JLSV}_1$) algorithm and ranking polynomials by their
3d Murphy E-score, after about 100 core hours found the following polynomial pair from the cyclic family
of degree six described in \cite{MR866107}:
\begin{flalign*}
f_0 &= x^6 - 40226000394x^5 - 100565001000x^4 - 20x^3 + 100565000985x^2\\
 &+ 40226000400x + 1 \\
f_1 &= 80447172120x^6 + 104483881186x^5 - 945497878835x^4 - 1608943442400x^3\\
 &- 261209702965x^2 + 378199151534x + 80447172120
\end{flalign*}
We computed the 3d alpha score for these and found $\alpha(f_0) = -3.6$ and $\alpha(f_1) = -12.6$.  

\subsection{Relation collection}

Our implementation was written as a standalone executable, independent of CADO-NFS, 
producing relations in the format that CADO-NFS can use.   We carry out cofactorization using
Pollard's $p-1$ algorithm and two rounds of Edwards elliptic curve factorization.  Although the
cofactorisation implementation uses a standard approach and is not an improvement on CADO-NFS's
cofactorization rig, our program is extremely fast to sieve.  This allowed us not only to use 
a larger factor base, but also to search for relations that are `twice as difficult' to find, i.e. to
use a large prime bound of $2^{28}$ as opposed to the $2^{29}$ used in the 422-bit record.
We were able to use a factor base bound of $2^{24}$ with no major loss of speed.  

In addition,
we were able to use a larger sieve region due to the speed of lattice enumeration.  We
used a sieving region of size $2^9 \times 2^{10} \times 2^{10}$, compared to the region
$2^{10} \times 2^{10} \times 2^8$ used in the previous record. 
The time per special-$\mathfrak{q}$ was roughly constant across the entire range, at between 150-170
seconds.  The bottleneck was cofactorization, by a wide margin - typically sieving takes less than
40 seconds per special-$\mathfrak{q}$, including CPU sorting (the list on each side typically has
about 400M-500M elements, each element taking 5 bytes.  Note that we omit the smallest primes
in the sieve as they correspond to dense lattices.  This is alleviated in trial factorization).
Cofactorization typically takes between 120-130 seconds.
We sieved most special-$\mathfrak{q}$s on the $f_0$ side with norm between 16M and 263M. We were
able to utilize all 40 cores on our sieving nodes, where each node has 192Gb of memory.  Our
program sieves only one ideal in each Galois orbit.  We apply the Galois automorphism as a 
post-processing step.  We found 7,152,855 unique relations and then applied the Galois automorphism
(which is trivial in core-hours) and after removing duplicates we were left with 34,115,391
unique relations.  The total sieving effort was 69,120 core hours.

\begin{table}[h]
\begin{center}
\begin{tabulary}{12cm}{|l|C|C|}
\hline
Authors&GGMT&This work\\
\hline
Field size (bits)&422&423\\
$\alpha$-values&-2.4,-14.3&-3.6,-12.6\\
Murphy-E&$2^{-20.51259}$&$2^{-20.45961}$\\
Sieving region $H$&10,10,8&9,10,10\\
Factor base bounds&$2^{21},2^{21}$&$2^{24},2^{24}$\\
Smoothness bounds&$2^{29},2^{29}$&$2^{28},2^{28}$\\
$\#S = q_{max}2^{H_0+H_1+H_2}$&$2^{55}$&$2^{57}$\\
Special-$\mathfrak{q}$ side&0&0\\
q-range&$]2^{21},2^{27.9}[$&$]2^{23.9},2^{30}[$\\
Galois action&6&6\\
\#unique relations&71,850,465&34,115,391\\
\#required relations&$\approx 56M$&29,246,136\\
purged&18,335,401&7,598,223\\
filtered&5,218,599&2,754,009\\
Total sieving time&201,600&69,120\\
\hline
\end{tabulary}
\end{center}
\caption{Key statistics of record computations in $\F_{p^6}$.  All timings in core hours}\label{RecordCompare}
\end{table}

\subsection{Construction of matrix}
We modified CADO-NFS \cite{CADO} to produce a matrix arising from degree-2 sieving ideals for the linear
algebra step.

\subsection{Linear algebra}
We used the Block Wiedemann implementation in CADO-NFS (we compiled commit \texttt{d6962f667d3c...}
with MPI enabled), with parameters $n = 10$ and $m = 20$.  Due to time constraints, we needed to 
minimize wall clock time so we chose to run the computation on 4 nodes, to reduce the iteration
time for the Krylov sequences.  Also, to avoid complications, we did not run the 10 Krylov
sequences in parallel.  The net result was that we spent 11,760 core hours on the Krylov step, which
is suboptimal (but got us the result in time). 
It took 24 core hours (on one core) to compute the linear generator and 672 core hours for the
solution step.  This gave $2,754,009$ of the factor base ideal virtual logarithms.
We ran the log reconstruction to give a final total of $25,215,976$ known
virtual logarithms out of a possible total of $29,246,136$ factor base ideals.\\

We note that a similar-sized matrix was solved in \cite{MR3775580}, which used 1,920 core hours for the Krylov step.
However, due to our choice of the large prime bound, set to $2^{28}$, our linear algebra effort
to set the new record was considerably less than that of the previous record of 422 bits,
which involved a Krylov step taking 23,390 core hours for a large prime bound set to $2^{29}$.

\subsection{Individual logarithm}
Take the element $g = x + 2 \in \F_{p^6} = \F_p[x]/\langle f_0(x) \rangle$. Let
\begin{flalign*}
\ell = 9589868090658955488259764600093934829209,
\end{flalign*}
a large prime factor of $p^2-p+1$. Let $h = \left(p^6-1\right)/\ell$.  Note that $g$ is not a
generator of the entire multiplicative subgroup of $\F_{p^6}$, but we do have that $g^h$ is a
generator of the subgroup of size $\ell$.  It is easy to compute $\text{vlog}(g)$ since
$N_0(g) = -3^3$.\\

\noindent We have $\text{vlog}(g) = 8951069617162908953536183274937613985265$.  We chose the target
{\small
\begin{flalign*}
t &= 314159265358979323846x^5 + 264338327950288419716x^4 + 939937510582097494459x^3 \\
  &+ 230781640628620899862x^2 + 803482534211706798214x + 808651328230664709384
\end{flalign*}
}%
We implemented the initial splitting algorithm of A. Guillevic \cite{MR3904147} in SAGE, and after a few core
hours found that
\begin{flalign*}
g^{74265}t = uvw(-129592286880919x^2 - 103570474976165x - 5550010113050)
\end{flalign*}
where $u \in \F_{p^2}, v \in \F_{p^3}, w \in \F_p$, so that their logarithm modulo $\ell$ is zero.
The norm of the latter term is $-11 \cdot 37 \cdot 71 \cdot 97 \cdot 197 \cdot 821 \cdot 24682829 \cdot$
$33769709 \cdot 83609989 \cdot 13978298429383 \cdot 21662603713879 \cdot 74293619085767 \cdot$
$141762919001833 \cdot 381566853770521$.  We had 5 special-$\mathfrak{q}$ to descend, the largest
having 49 bits.  We used our 3d lattice sieve implementation to descend from these ideals of unknown log to factor
base elements with known logarithms.  This was a somewhat manual process and took about a 
day's work (about 8 hours).\\

\noindent We obtained $\text{vlog}(t) = 2619623637064116359346428467068287245870$, so that
\begin{flalign*}
\text{log}_g(t) \equiv \text{vlog}(t)/\text{vlog}(g) \equiv 7435826750517015269718230402645557947880\, \text{mod}\, \ell.
\end{flalign*}

\begin{table}[h]
\begin{center}
\begin{tabulary}{12cm}{|l|C|l|l|C|C|C|}
\hline
year&size of $p^n$&authors&algorithm&rel. col.&lin. alg.&total\\
\hline
2008&240&Zajac&NFS-HD&580&322&912\\
2015&240&HAKT&NFS-HD&527&-&-\\
2017&240&GGMT&NFS-HD&22&5&27\\
2017&300&GGMT&NFS-HD&164&39&203\\
2017&389&GGMT&NFS-HD&18,960&2,400&21,360\\
2017&422&GGMT&NFS-HD&201,600&26,880&228,480\\
2019&423&this work&NFS-HD&69,120&12,480&81,600\\
\hline
\end{tabulary}
\end{center}
\caption{Comparison with other record computations in $\F_{p^6}$.  All timings in core hours}\label{RecordCompare2}
\end{table}

\section{Pairing break}

Let $p$ be the same prime as in the previous section.
Define $\F_{p^2} = \F_p[i]/\langle{i^2+2\rangle}$.  The curve $E/\F_{p^2} : y^2 = x^3 + b$, $b = i + 7$ is
supersingular of trace $p$, hence of order $p^2-p+1$. Define $\F_{p^6} = \F_{p^2}[j]/\langle{j^3-b \rangle}$.
The embedding field of the curve $E$ is $\F_{p^6}$.  We take
$$G_0 = \left(5, 751568328314480688740i +751642554083315688493\right)$$
 and we check that $G_1 = [273]G_0$ is a generator of $E(\F_{p^2})[\ell]$.  The distortion map
$\phi : (x,y) \mapsto \left(x^p / \left(jb^{(p-2)/3}\right), y^p/\left(b^{(p-1)/2}\right)\right)$ gives
a generator $G_2 = \phi(G_1)$ of the second dimension of the $\ell$-torsion.  We take the point
$$P_0 = (314159265358979323846i + 264338327950288419717, $$
$$ 1560320966141767888064i + 368067364535991558380)$$
from the decimals of $\pi$, and $P_1 = [273]P_0 \in E(\F_{p^2})[\ell]$ is our challenge.  We aim to
compute the discrete logarithm of $P_1$ to base $G_1$.  To do so, we transfer $G_1$ and $P_1$
to $\F_{p^6}$, and obtain $g = e_{Tate}\left(G_1, \phi(G_1)\right)$ and $t = e_{Tate}\left(P_1,\phi(G_1)\right)$, or
{\small
\begin{flalign*}
t &= 709659446396572245219x^5 + 760855550263311226560x^4 + 459517758627469463106x^3\\
 &+ 1075867962756498791880x^2 + 966415406496231787507x + 759380554536416832249,\\
g &= 1445115464416256318145x^5 + 608219705720308630653x^4 + 1328213831161031326049x^3\\
 &+ 104723931403852502861x^2 + 1118264722333528462011x + 551285267384030855316
\end{flalign*}
}%
The initial splitting gave a 50-bit smooth generator
\begin{flalign*}
g^{289236} = uvw\left(-207659249318101x^2 - 32084626907475x + 36052674649889\right)
\end{flalign*}
where $u \in \F_{p^2}, v \in \F_{p^3}, w \in \F_p$, so that their logarithm modulo $\ell$ is zero.
The norm of the latter term is $11 \cdot 71 \cdot 79 \cdot 1453 \cdot 433123 \cdot 85478849 \cdot 34588617703 \cdot$
{\small$ 40197196124443 \cdot 76694584420127 \cdot 370667620290007 \cdot 419573910884273 \cdot 823157513981483$}.
We had 6 special-$\mathfrak{q}$ to descend.
We also got a 49-bit smooth challenge of norm $23 \cdot 29^2 \cdot 41 \cdot 563 \cdot 2917 \cdot 1245103 \cdot 12006859 \cdot$
$107347203833 \cdot 506649149393 \cdot 39018481981309 \cdot 138780153403907 \cdot 174514280440993 \cdot 302260510161053$:
\begin{flalign*}
g^{91260}t = uvw\left(-59788863574984x^2 + 62066870577408x + 88384197770333\right)
\end{flalign*}
We obtained $\text{vlog}(g) = 7599151482912535295281621925658364195913$ and\\
$\text{vlog}(t) = 4642225023760573112152590887355819325364$, so that
$\text{log}_g(t) \equiv \text{vlog}(t)/\text{vlog}(g) \equiv 4325953856049730257332335443497115431763\,\text{mod}\,\ell$.


\section{Conclusion}
We have presented a new approach to lattice sieving in higher dimensions for the number field sieve,
together with a novel approach to avoiding inefficient memory access patterns which applies
regardless of dimension.  In addition, we implemented the 3d case of our idea and used it to
set a record in solving discrete log in $\F_{p^6}$, a typical target in cryptanalysis of
pairing-based cryptography, in time a factor of more than $2.5\times$ better (in core hours) 
than the previous record, which was of a directly comparable size.
It should be possible to improve the code further with more
effort put into optimization.  We have indicated that the sieving enumeration generalizes to
higher dimensions as long as a certain integer linear programming problem is tractable.  This
has immediate implications for the possibility of implementation of the Tower Number
Field Sieve and e.g. the Extended Tower Number Field Sieve, the latter of which is dependent
on sieving in dimension at least four.  The recent preprint \cite{guillevic:hal-02263098}
addresses one major prerequisite to the realization of TNFS and ExTNFS, concerning polynomial
selection, while in the present work we give a strong indication that another obstruction,
that of sieving efficiently in small dimensions of four and above, may be easier than
 first thought.

\newpage

\nocite{*}
\bibliographystyle{splncs04}
\bibliography{refs}

\end{document}